\documentclass[10pt]{amsart}

\usepackage[T1]{fontenc}
\usepackage{amsmath,amssymb,amsthm}
\usepackage{microtype}
\usepackage[margin=0.6in]{geometry}

\newtheorem{theorem}{Theorem}[section]
\newtheorem{corollary}[theorem]{Corollary}

\title[Rank-One Sutured Instanton Homology]{Rank-One Sutured Instanton Homology Detects Homology Products}
\author{Qilong Guo}
\address{Department of Mathematics, College of Science, China University of Petroleum-Beijing, Beijing 102249, China}
\email{guoqilong@cup.edu.cn}
\date{}
\subjclass[2020]{Primary 57K31; Secondary 57R58}
\keywords{Sutured instanton homology, sutured torsion, product sutured manifold, homology product}

\begin{document}

\begin{abstract}
We prove that a balanced sutured manifold \((M,\gamma)\) satisfying
\(\dim_{\mathbb C} SHI(M,\gamma)=1\) is an integral homology product; that is,
\[
H_*(M,R_-(\gamma);\mathbb Z)
=
H_*(M,R_+(\gamma);\mathbb Z)=0.
\]
Consequently, for every taut balanced sutured manifold \((M,\gamma)\),
\[
SHI(M,\gamma)\cong\mathbb C
\quad\Longleftrightarrow\quad
(M,\gamma)\text{ is a product sutured manifold}.
\]
The proof uses the enhanced Euler characteristic of \(SHI\) and a
specialization formula for sutured torsion.
\end{abstract}

\maketitle

\section{Rank One and Homology Products}

For a sutured manifold \((M,\gamma)\), let \(SHI(M,\gamma)\) denote
its sutured instanton homology. Kronheimer and Mrowka proved in
\cite[Theorem~7.18]{KM} that a taut balanced sutured manifold
\((M,\gamma)\) which is a homology product is a product sutured manifold
if and only if
\[
SHI(M,\gamma)\cong\mathbb C.
\]
We show that the homology-product hypothesis already follows from the
rank-one assumption. The key input is the enhanced Euler characteristic
constructed by Li and Ye \cite{LY}.

Throughout, homology groups have integral coefficients unless indicated
otherwise. We call \((M,\gamma)\) an integral homology product if the
inclusions
\[
R_\pm(\gamma)\hookrightarrow M
\]
induce isomorphisms on integral homology. Equivalently,
\[
H_*(M,R_-(\gamma);\mathbb Z)
=
H_*(M,R_+(\gamma);\mathbb Z)
=0.
\]
\begin{theorem}
Let \((M,\gamma)\) be a balanced sutured manifold. If
\[
\dim_{\mathbb C}SHI(M,\gamma)=1,
\]
then \((M,\gamma)\) is an integral homology product.
\end{theorem}

\begin{proof}
Write \(H=H_1(M;\mathbb Z)\). Applying \cite[Theorem~5.1]{LY} to
\((-M,-\gamma)\) gives
\begin{equation}
\chi_{\mathrm{en}}\bigl(SHI(M,\gamma)\bigr)
=
\tau(M,\gamma)
\quad\text{in }\mathbb Z[H]/\!\pm H.
\tag{1}
\end{equation}
Here \(\chi_{\mathrm{en}}\) denotes the enhanced Euler characteristic,
\(\tau(M,\gamma)\) denotes the sutured torsion, and \(\mathbb Z[H]\)
denotes the integral group ring of \(H\). The notation
\(\mathbb Z[H]/\!\pm H\) identifies elements that differ by
multiplication by a unit \(\pm h\), where \(h\in H\).

In the \(H\)-indexed decomposition used to define
\(\chi_{\mathrm{en}}\), there is a unique nonzero summand, which is
one-dimensional. Hence \(\chi_{\mathrm{en}}(SHI(M,\gamma))\) has a
representative that is a monomial with coefficient \(\pm1\) in
\(\mathbb Z[H]\). By (1), the class of \(\tau(M,\gamma)\) has such a
representative. Thus, there exists a representative
\(\widetilde{\tau}(M,\gamma)\in\mathbb Z[H]\) of the class
\(\tau(M,\gamma)\) such that
\begin{equation}
\widetilde{\tau}(M,\gamma)=\pm h_0
\quad\text{for some }h_0\in H.
\tag{2}
\end{equation}

Let
\[
p:H_1(M;\mathbb Z)\longrightarrow
G:=H_1(M,R_-(\gamma);\mathbb Z)
\]
be the natural map. Write
\[
p_*:\mathbb Z[H]\longrightarrow\mathbb Z[G]
\]
for the induced map of group rings. Proposition~6.5 of \cite{FJR} gives,
up to multiplication by an element of \(\pm G\),
\begin{equation}
p_*\bigl(\widetilde{\tau}(M,\gamma)\bigr)=\pm I_G,
\qquad
I_G=
\begin{cases}
\displaystyle\sum_{g\in G}g, & G\text{ is finite},\\[4pt]
0_{\mathbb Z[G]}, & G\text{ is infinite}.
\end{cases}
\tag{3}
\end{equation}
Here \(\mathbb Z[G]\) is the integral group ring of \(G\), and \(I_G\)
is the indicated element of that group ring.

By (2), the left-hand side of (3) is the nonzero monomial
\(\pm p(h_0)\). Hence \(G\) is finite. The support of \(I_G\) then has
cardinality \(|G|\), so (3) forces \(|G|=1\). Thus,
\begin{equation}
H_1(M,R_-(\gamma);\mathbb Z)=0.
\tag{4}
\end{equation}

It remains to check the other relative homology groups. Since the
annular and toroidal components of \(\gamma\) have Euler characteristic
zero, we have
\[
2\chi(M)=\chi(\partial M)
=\chi(R_+(\gamma))+\chi(R_-(\gamma)),
\]
where \(\chi\) denotes the Euler characteristic. Since \((M,\gamma)\)
is balanced, we have
\[
\chi(R_+(\gamma))=\chi(R_-(\gamma)).
\]
It follows that \(\chi(M)=\chi(R_-(\gamma))\). By the definition of
relative Euler characteristic,
\[
\chi(M,R_-(\gamma))
=\chi(M)-\chi(R_-(\gamma))=0,
\]
where \(\chi(M,R_-(\gamma))\) denotes the relative Euler characteristic
of the pair \((M,R_-(\gamma))\). The balanced condition also implies that every component of
\(M\) meets both \(R_-(\gamma)\) and \(R_+(\gamma)\). Hence the
inclusion-induced maps
\[
H_0(R_\pm(\gamma);\mathbb Z)\longrightarrow H_0(M;\mathbb Z)
\]
are surjective, and therefore
\[
H_0(M,R_\pm(\gamma);\mathbb Z)=0.
\]
Poincar\'e--Lefschetz duality gives
\[
H_3(M,R_-(\gamma);\mathbb Z)
\cong H^0(M,R_+(\gamma);\mathbb Z)=0.
\]
Together with (4), the relative Euler characteristic shows that the
second relative homology group has rank zero. On the other hand,
duality and the universal coefficient theorem give
\[
\begin{aligned}
H_2(M,R_-(\gamma);\mathbb Z)
&\cong H^1(M,R_+(\gamma);\mathbb Z)\\
&\cong
\operatorname{Hom}\bigl(
H_1(M,R_+(\gamma);\mathbb Z),\mathbb Z
\bigr).
\end{aligned}
\]
Since \(H_1(M,R_+(\gamma);\mathbb Z)\) is finitely generated,
the last group is free abelian. Since the first group has rank zero,
all three groups vanish.
Thus
\[
H_*(M,R_-(\gamma);\mathbb Z)=0.
\]
Duality now gives
\[
H^*(M,R_+(\gamma);\mathbb Z)=0,
\]
and, since the relative homology groups are finitely generated,
the universal coefficient theorem, applied degree by degree, gives
\[
H_*(M,R_+(\gamma);\mathbb Z)=0
\]
as well. Hence \((M,\gamma)\) is an integral homology product.
\end{proof}

\begin{corollary}
Let \((M,\gamma)\) be a taut balanced sutured manifold. Then
\[
SHI(M,\gamma)\cong\mathbb C
\quad\Longleftrightarrow\quad
(M,\gamma)\text{ is a product sutured manifold}.
\]
\end{corollary}

\begin{proof}
The forward implication follows from Theorem 1.1 and \cite[Theorem~7.18]{KM};
the reverse implication is the remark following \cite[Definition~7.10]{KM}.
\end{proof}

\section*{Acknowledgements}
The author thanks Zhenkun Li for helpful discussions.

\section*{Use of AI Tools}
AI tools were used solely for language polishing. The author takes
full responsibility for the mathematical content of this paper.

\normalsize


\begin{thebibliography}{9}
\small

\bibitem{FJR}
S.~Friedl, A.~Juh\'asz, and J.~Rasmussen,
\newblock The decategorification of sutured Floer homology,
\newblock \emph{J. Topol.} \textbf{4} (2011), no.~2, 431--478.

\bibitem{KM}
P.~B. Kronheimer and T.~S. Mrowka,
\newblock Knots, sutures, and excision,
\newblock \emph{J. Differential Geom.} \textbf{84} (2010), no.~2, 301--364.

\bibitem{LY}
Z.~Li and F.~Ye,
\newblock An enhanced Euler characteristic of sutured instanton homology,
\newblock \emph{Int. Math. Res. Not. IMRN} 2024 (2024), no.~4, 2873--2936.

\end{thebibliography}
\end{document}